# Maximization of the Spectral Gap for Chemical Graphs by means of a Solution to a Mixed Integer Semidefinite Program

*Soňa Pavlíková and Daniel Ševčovič\**

*FCFT, Slovak Technical University, 812 37 Bratislava, Slovakia, and
FMFI, Comenius University, 842 48 Bratislava, Slovakia,*
*\*Corresponding author: sevcovic@fmph.uniba.sk*

**Abstract.** In this paper we analyze the spectral gap of a weighted graph which is the difference between the smallest positive and largest negative eigenvalue of its adjacency matrix. Such a graph can represent e.g. a chemical organic molecule. Our goal is to construct a new graph by bridging two given weighted graphs over a bipartite graph. The aim is to maximize the spectral gap with respect to a bridging graph. To this end, we construct a mixed integer semidefinite program for maximization of the spectral gap and compute it numerically.

## 1. Introduction

Following (Streitwieser, 1961), eigenvalues of a graph describing an organic molecule are related to energies of molecular orbitals (see also Li et al., 2013; Zhang & An, 1999; (Fowler et al., 2001). In computational chemistry, Hückel's molecular orbital method (Hückel, 1931) is applied in order to analyze stability of chemical molecules (Streitwieser, 1961), (Gutman & Rouvray, 1979). According to this method, the energies $E_k, k=1,\cdots,n,$ of an organic molecule are eigenvalues of the Hamiltonian matrix $H$ and eigenvectors of $H$ are called orbitals. The structure of the square symmetric matrix $H$ is as follows: $H_{ii} = \alpha$ for the carbon C atom located at the $i$-th vertex, and $H_{ii} = \alpha + h_A \beta$ for other atoms A, where $\alpha < 0$ is the Coulomb integral and $\beta < 0$ is the resonance integral, $H_{ij} = \beta$ if both vertices $i$ and $j$ are carbon C atoms, $H_{ij} = k_{AB} \beta$ for other neighboring atoms A and B, $H_{ij} = 0$ otherwise. The atomic constants $h_A$ and $k_{AB}$ should be specified. In the case of a hydrocarbon we have $h_C = 0$ and $k_{CC} = 1$. On the other hand, the molecule of pyridine contains one atom of nitrate N and five atoms of carbon C (see figure 1, c)). Nitrate N has the higher Pauling electronegativity with respect to the carbon C atom. With regard to (Streitwieser, 1961, Chapter 5.1), $h_N = 0.5 > 0$ and $k_{CN} = 0.8$. Clearly, in the case of a hydrocarbon molecule we have $H = \alpha I + \beta A$ where $I$ is the identity matrix and $A$ is the binary adjacency matrix of the molecular structural graph $G$. Hence $E_k = \alpha + \beta \lambda_k$ where $\lambda_k$ is an eigenvalue of $A$.

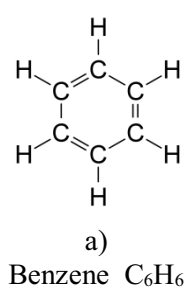 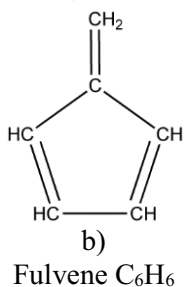 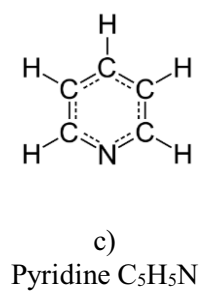

a)  
Benzene C$_6$H$_6$

b)  
Fulvene C$_6$H$_6$

c)  
Pyridine C$_5$H$_5$N

**Figure 1** Chemical graphs of benzene, fulvene and pyridine.

In three examples of organic molecules shown in figure 1, the Hamiltonian matrix has the form $H = \alpha I + \beta A$ where the matrix $A$ is defined as follows:

$$A_{Benzene} = \begin{pmatrix} 0 & 1 & 0 & 0 & 0 & 1 \\ 1 & 0 & 1 & 0 & 0 & 0 \\ 0 & 1 & 0 & 1 & 0 & 0 \\ 0 & 0 & 1 & 0 & 1 & 0 \\ 0 & 0 & 0 & 1 & 0 & 1 \\ 1 & 0 & 0 & 0 & 1 & 0 \end{pmatrix}, \quad A_{Fulvene} = \begin{pmatrix} 0 & 1 & 0 & 0 & 1 & 0 \\ 1 & 0 & 1 & 0 & 0 & 0 \\ 0 & 1 & 0 & 1 & 0 & 0 \\ 0 & 0 & 1 & 0 & 1 & 1 \\ 1 & 0 & 0 & 1 & 0 & 0 \\ 0 & 0 & 0 & 1 & 0 & 0 \end{pmatrix}, \quad A_{Pyridine} = D^A \begin{pmatrix} 1 & 1 & 0 & 0 & 0 & 1 \\ 1 & 0 & 1 & 0 & 0 & 0 \\ 0 & 1 & 0 & 1 & 0 & 0 \\ 0 & 0 & 1 & 0 & 1 & 0 \\ 0 & 0 & 0 & 1 & 0 & 1 \\ 1 & 0 & 0 & 0 & 1 & 0 \end{pmatrix} D^A.$$

Here $D^A = diag(\sqrt{h_N}, \frac{k_{CN}}{\sqrt{h_N}}, 1, 1, 1, \frac{k_{CN}}{\sqrt{h_N}})$ is the diagonal matrix.

In what follows, a non-oriented weighted graph $G_A$ on $n$ vertices with a weighted symmetric adjacency matrix $A$ is called a *voltage graph* if there exists a diagonal matrix $D^A$ having nonzero elements and such that $A = D^A \overline{A} D^A$ where $\overline{A}$ is a binary adjacency matrix of the graph $G_A$, i.e. $\overline{A}_{ij} = 1$ if there exists an edge from the vertex $i$ to the vertex $j$, $\overline{A}_{ij} = 0$ otherwise. The spectrum $\sigma(G)$ of an undirected weighted graph $G$ consists of eigenvalues of its $n \times n$ symmetric weighted adjacency matrix $A(G)$, i.e. $\sigma(G) = \{\lambda_k, k = 1, \cdots, n\}$ where $\lambda_1 \geq \cdots \geq \lambda_n$ are eigenvalues of $A(G)$, (c.f. Brouwer & Haemers, 2012; Cvetković et al., 1980). If the spectrum does not contain zero then there exists the inverse matrix $A^{-1}$ of the adjacency matrix $A = A(G)$. According to (Godsil, 1985) a graph $G_A$ is called invertible if the inverse matrix $A^{-1}$ is diagonally similar to a nonnegative matrix representing the inverse graph $G^{-1}$. In computational chemistry (c.f. Streitwaiser, 1961), the energy $E_{HOMO}$ of the highest occupied molecular orbital (HOMO) corresponds to the eigenvalue $\lambda_{HOMO} = \lambda_k$ where $k = n/2$ for $n$ even and $k = (n+1)/2$ for $n$ odd. The energy $E_{LUMO}$ of the lowest unoccupied molecular orbital (LUMO) corresponds to the subsequent eigenvalue $\lambda_{LUMO} = \lambda_{k+1}$ for $n$ even, and $\lambda_{LUMO} = \lambda_k$ for $n$ odd. The HOMO-LUMO separation gap is the difference between $E_{LUMO}$ and $E_{HOMO}$ energies, i.e. $E_{LUMO} - E_{HOMO} = -\beta(\lambda_{HOMO} - \lambda_{LUMO}) \geq 0$ because $\beta < 0$. For molecules with closed shells (Fowler & Pisański, 2010) we have $\lambda_{HOMO} > 0 > \lambda_{LUMO}$. In this case, the HOMO-LUMO separation gap is equal to the energy difference $E_{LUMO} - E_{HOMO} = -\beta \Lambda_{HL}(G_A)$ where

$$\Lambda_{HL}(G_A) = \lambda^+(G_A) - \lambda^-(G_A). \tag{1}$$

Here $\lambda^+(G_A)$ is the smallest positive eigenvalue and $\lambda^-(G_A)$ is the largest negative eigenvalue of the adjacency matrix $A$ of the structural molecular graph $G_A$. Following (Aihara, 1999a,b), it is energetically unfavorable to add electrons to a high-lying LUMO orbital. Hence a larger HOMO-LUMO gap implies higher kinetic stability and low chemical reactivity of a molecule. According to (Bacalis & Zdetsis, 2009) the HOMO-LUMO energy gap is generally decreasing with the number $n$ of vertices of the structural graph.

In this paper, our goal is to investigate extremal properties of the HOMO-LUMO spectral gap $\Lambda_{HL}(G_A)$. We show how to represent $\Lambda_{HL}(G_A)$ by means of an optimal solution to a mixed integer semidefinite programming problem. Our aim is to investigate dependence of the HOMO-LUMO spectral gap for a class of structural voltage graphs which can be constructed from two given voltage graphs by bridging them over a bipartite voltage graph.

## 2. Semidefinite programming representation of the HOMO-LUMO spectral gap

Clearly, the smallest positive and largest negative eigenvalues of $G_C$ can be expressed as follows:

$$\lambda^+(G_C) = \frac{1}{\lambda_{max}(C^{-1})}, \qquad \lambda^-(G_C) = \frac{1}{\lambda_{min}(C^{-1})}, \tag{2}$$

where $\lambda_{max}(C^{-1}) > 0$ and $\lambda_{min}(C^{-1}) = -\lambda_{max}(-C^{-1}) < 0$ are the maximum and minimum eigenvalues of the inverse matrix $C^{-1}$, respectively. In what follows, we denote by $\leq$ the Löwner partial ordering on symmetric matrices, i.e. $A \leq B$ if and only if the matrix $B - A$ is positive semidefinite, that is $B - A \geq 0$. Then the maximal and minimal eigenvalues of $C^{-1}$ can be expressed as follows:

$$0 < \lambda_{max}(C^{-1}) = \min_{C^{-1} \leq tI} t, \qquad 0 > \lambda_{min}(C^{-1}) = \max_{sI \leq C^{-1}} s, \tag{3}$$

(see Boyd & Vanderberghe, 2004; Cvetković et al. 2004). Applying the change of variables: $\mu = 1/t, \eta = -1/s$, we obtain the following characterization of the lowest positive and largest negative eigenvalues of a graph $G_C$:

$$\lambda^+(G_C) = \max_{\mu C^{-1} \leq I} \mu, \qquad \lambda^-(G_C) = -\max_{-\eta C^{-1} \leq I} \eta. \tag{4}$$

As a consequence, we obtain the following semidefinite representation of the HOMO-LUMO spectral gap. Let $G_C$ be a vertex labeled graph. Then the HOMO-LUMO spectral gap $\Lambda_{HL}(G_C)$ of $G_C$ is the optimal value of the following semidefinite programming problem:

$$\begin{aligned}\Lambda_{HL}(G_C) = \max_{\mu, \eta \geq 0} \quad & \mu + \eta \\ s.t. \quad & \mu C^{-1} \leq I, \quad -\eta C^{-1} \leq I.\end{aligned} \tag{5}$$

## 3. Bridged graphs and their HOMO-LUMO spectral gap

In this section we recall a notion of a voltage graph which is constructed from two given voltage graphs $G_A$ and $G_B$ by bridging vertices of $G_A$ to vertices of $G_B$ through a bipartite graph $G_H$ with the adjacency matrix:

$$A(G_H) = \begin{pmatrix} 0 & H \\ H^T & 0 \end{pmatrix}, \tag{6}$$

where $H$ is an $n \times m$ matrix. More precisely, let $G_A$ and $G_B$ be two undirected vertex-labeled voltage graphs on $n$ and $m$ vertices, respectively. By $B_H(G_A, G_B)$ we shall denote the graph $G_C$ on $n+m$ vertices which is obtained by bridging vertices $\{1, \cdots, n\}$ of the vertex-labeled graph $G_A$ to vertices $\{1, \cdots, m\}$ of $G_B$ through an $(n,m)$-bipartite graph $G_H$. It means that the adjacency matrix $C = A(G_C)$ of $G_C$ has the form:

$$C = \begin{pmatrix} A & H \\ H^T & B \end{pmatrix} \qquad (7)$$

(c.f. Pavlíková & Ševčovič, D., 2016a,b). In what follows, we will assume that the graphs $G_A$ and $G_B$ are voltage graphs, i.e. there exist diagonal matrices $D^A$ and $D^B$ such that $A = D^A \overline{A} D^A$ and $B = D^B \overline{B} D^B$ where $\overline{A}$ and $\overline{B}$ are binary adjacency matrices. Concerning the $n \times m$ matrix $H$ we shall assume that $H = D^A \overline{H} D^B$ where $\overline{H}$ is a binary matrix. Moreover, matrices $A$ and $B$ are assumed to be invertible. Regarding the graph $G_B$ we shall furthermore assume that $G_B$ is arbitrarily bridgeable over the first $\{1, \cdots, k_B\}$ vertices of $G_B$. It means that the $k_B \times k_B$ upper principal sub-matrix is a null matrix, i.e. $(B^{-1})_{ij} = 0$ for all $i, j = 1, \cdots, k_B$ (c.f. Pavlíková, 2016, Pavlíková & Ševčovič, D., 2016a,b). For example, a graph representing the fulvene molecule (see figure 1, b)) has the inverse matrix:

$$A^{-1}_{Fulvene} = \begin{pmatrix} \boxed{0 & 0} & 0 & 0 & 1 & -1 \\ 0 & 0 & 1 & 0 & 0 & -1 \\ 0 & 1 & 0 & 0 & -1 & 1 \\ 0 & 0 & 0 & 0 & 0 & 1 \\ 1 & 0 & -1 & 0 & 0 & 1 \\ -1 & -1 & 1 & 1 & 1 & -2 \end{pmatrix},$$

which means that it is arbitrarily bridgeable over the vertices $\{1, 2\}$. Similarly, a graph representing the organic molecule of benzene (see figure 1, a)) is arbitrarily bridgeable (after suitable permutation of its vertices) over any pair of vertices having the distance 2.

Now, with regard to the Schur complement theorem, we obtain:

$$C^{-1} = \begin{pmatrix} A & H \\ H^T & B \end{pmatrix}^{-1} = \begin{pmatrix} S^{-1} & -S^{-1}HB^{-1} \\ -B^{-1}H^T S^{-1} & B^{-1} + B^{-1}H^T S^{-1} H B^{-1} \end{pmatrix} = Q^T \begin{pmatrix} S^{-1} & 0 \\ 0 & B^{-1} \end{pmatrix} Q, \qquad (8)$$

where $S = A - HB^{-1}H^T$ is the Schur complement and $Q$ is the invertible matrix with $Z = Q^{-1}$ given by:

$$Q = \begin{pmatrix} I & -HB^{-1} \\ 0 & I \end{pmatrix}, \qquad Z = \begin{pmatrix} I & HB^{-1} \\ 0 & I \end{pmatrix}, \qquad (9)$$

(see Hamala & Trnovská, 2013). Now, let $G_C = B_H(G_A, G_B)$ be the graph obtained from graphs $G_A$ and $G_B$ by bridging them through a bipartite graph $G_H$ with the adjacency matrix $A(G_H)$ having the off-diagonal block sub-matrix $H$ (see (7)) with the property:

$$H_{ij} = 0 \text{ for all } i = 1, \cdots, n \text{ and } j = k_B + 1, \cdots, m. \qquad (10)$$

Since the graph $G_B$ is assumed to be arbitrarily bridgeable over the first $k_B$ vertices and the matrix $H$ has the property (10) then $HB^{-1}H^T = 0$. Therefore, for the Schur complement we have $S = A - HB^{-1}H^T = A$. Hence $S^{-1} = A^{-1}$. The inverse matrix $C^{-1}$ of the adjacency matrix $C$ is given by:

$$C^{-1} = \begin{pmatrix} A^{-1} & -A^{-1}HB^{-1} \\ -B^{-1}H^T A^{-1} & B^{-1}+B^{-1}H^T A^{-1}HB^{-1} \end{pmatrix} = Q^T \begin{pmatrix} A^{-1} & 0 \\ 0 & B^{-1} \end{pmatrix} Q.$$

Then, for $\mu \geq 0$, we have $\mu C^{-1} \leq I$ if and only if $\mu Z^T C^{-1} Z \leq Z^T Z$, i.e.,

$$\mu \begin{pmatrix} A^{-1} & 0 \\ 0 & B^{-1} \end{pmatrix} \leq Z^T Z = \begin{pmatrix} I & HB^{-1} \\ B^{-1}H^T & I+B^{-1}H^T HB^{-1} \end{pmatrix}.$$

This way we obtain the following representation of the HOMO-LUMO spectral gap $\Lambda_{HL}(G_C)$ for the bridged graph $G_C = B_H(G_A, G_B)$:

$$\Lambda_{HL}(G_C) = \max_{\mu,\eta \geq 0} \quad \mu + \eta$$

$$\begin{pmatrix} I - \mu A^{-1} & HB^{-1} \\ B^{-1}H^T & I - \mu B^{-1}+B^{-1}H^T HB^{-1} \end{pmatrix} \geq 0, \quad \begin{pmatrix} I + \eta A^{-1} & HB^{-1} \\ B^{-1}H^T & I + \eta B^{-1}+B^{-1}H^T HB^{-1} \end{pmatrix} \geq 0. \quad (11)$$

## 4. Maximization of the HOMO-LUMO spectral gap

Next we focus our attention on extremal properties of the HOMO-LUMO spectral gap for bridged graphs. Let $G_A$ be an invertible graph and $G_B$ be an arbitrarily bridgeable invertible graph. Our goal is to find a matrix $H$ forming the bipartite graph $G_H$ (see (6)) such that the HOMO-LUMO spectral gap $\Lambda_{HL}(G_C)$ is maximal, where $G_C = B_H(G_A, G_B)$. With regard to (11), the maximal HOMO-LUMO gap $\Lambda_{HL}^{opt} = \Lambda_{HL}^{opt}(G_A, G_B)$ with respect to a bipartite matrix $H$ can be obtained as the optimal value of the following mixed integer nonlinear optimization problem:

$$\Lambda_{HL}^{opt} = \max_{\substack{\mu,\eta \geq 0 \\ H,W}} \quad \mu + \eta$$

$$\begin{pmatrix} I - \mu A^{-1} & HB^{-1} \\ B^{-1}H^T & I - \mu B^{-1}+B^{-1}WB^{-1} \end{pmatrix} \geq 0, \quad \begin{pmatrix} I + \eta A^{-1} & HB^{-1} \\ B^{-1}H^T & I + \eta B^{-1}+B^{-1}WB^{-1} \end{pmatrix} \geq 0. \quad (12)$$

$W = H^T H, \quad H = D^A \overline{H} D^B, \quad \overline{H}_{ij} \in \{0,1\} \quad \text{for all} \quad i,j, \quad \sum_{k,l} \overline{H}_{kl} \geq 1.$

The last constraint $\sum_{k,l} \overline{H}_{kl} \geq 1$ means that there exists a bridge from $G_A$ to $G_B$. In the optimization problem (12), the objective function and the first two matrix inequality constraints are linear in the variables $\mu, \eta, H, W$. The constraint $W = H^T H$ makes the problem nonconvex. Moreover, there is a binary constraint $\overline{H}_{ij} \in \{0,1\}$. It means that (12) is a mixed integer nonconvex programming problem which is, in general, NP-hard to solve. In the last decades, various relaxation techniques for solving a mixed integer nonconvex problem have been developed (Boyd & Vanderberghe, 2004). Usually, semidefinite relaxations of an original nonconvex problem can be constructed by means of the second Lagrangian dual problem which is a convex semidefinite programming problem (c.f. Ševčovič & Trnovská, 2016). The binary constraint $\overline{H}_{ij} \in \{0,1\}$ is equivalent to the equality: $\overline{H}_{ij} = \overline{H}_{ij}^2$. Hence, for $W = H^T H$ where $H = D^A \overline{H} D^B, \overline{H}_{ij} \in \{0,1\}$, we have

$$W_{jj} = \sum_l H_{lj}^2 = \sum_l (D_{ll}^A)^2 \overline{H}_{lj}^2 (D_{jj}^B)^2 = \sum_l D_{ll}^A H_{lj} D_{jj}^B \quad \text{for all } j = 1, \cdots, m, \qquad (13)$$

which is a linear constraint between matrices $W$ and $H$. Now the nonconvex constraint $W = H^T H$ can be relaxed by a convex matrix inequality constraint: $W \geq H^T H$ which is suitable for application of available semidefinite programming solvers. Notice that the relaxation $W \geq H^T H$ together with the linear constraint (13) is tight in the following sense: if $W \geq H^T H$ and $W_{jj} = \sum_l D_{ll}^A H_{lj} D_{jj}^B$ for all $j = 1, \cdots, m$, where $H = D^A \overline{H} D^B, \overline{H}_{ij} \in \{0,1\}$, then $W = H^T H$. Indeed, let us denote. Then $L \geq 0$ and, with regard to (13), we have $L_{jj} = W_{jj} - \sum_l H_{lj}^2 = W_{jj} - \sum_l (D_{ll}^A)^2 \overline{H}_{lj}^2 (D_{jj}^B)^2 = W_{jj} - \sum_l D_{ll}^A H_{lj} D_{jj}^B = 0$ for all $j = 1, \cdots, m$, i.e. $diag(L) = 0$. It means that $L = 0$ and so $W = H^T H$. Hence the nonconvex-integer programming problem (12) can be solved by means of the following mixed integer semidefinite programming problem with linear matrix inequality constraints and integer constraints:

$$\Lambda_{HL}^{opt} = \max_{\substack{\mu, \eta \geq 0 \\ H, W}} \mu + \eta$$

$$\begin{pmatrix} I - \mu A^{-1} & HB^{-1} \\ B^{-1}H^T & I - \mu B^{-1} + B^{-1}WB^{-1} \end{pmatrix} \geq 0, \quad \begin{pmatrix} I + \eta A^{-1} & HB^{-1} \\ B^{-1}H^T & I + \eta B^{-1} + B^{-1}WB^{-1} \end{pmatrix} \geq 0, \qquad (14)$$

$W \geq H^T H, \quad H = D^A \overline{H} D^B, \quad \overline{H}_{ij} \in \{0,1\}, \quad W_{jj} = \sum_l D_{ll}^A H_{lj} D_{jj}^B$

for all $i, j$, and $\sum_{k,l} \overline{H}_{kl} \geq 1$.

In order to compute the mixed integer semidefinite program (14) one can use SeDuMi solver developed by (Sturm, 1999) within the Yalmip Matlab programming framework due to (Löfberg, 2004). The Matlab code for implementation of (14) is shown in table 1 below.

**Table 1.** The code for solving the mixed integer semidefinite programming problem (14).

```
% enter data:
% n x n matrix A, m x m matrix B, voltage matrices DA and DB, maxdegree

D=[[DA, zeros(n,m)]; [zeros(m,n),DB]];
% define data types
mu=sdpvar(1); eta=sdpvar(1); W=sdpvar(m,m);
H=sdpvar(n,m, 'full'); Haux=sdpvar(n,m, 'full'); Htilde=binvar(n,m, 'full');

% select Branch & Bound mixed integer programming solver for Yalmip
options=sdpsettings('solver','bnb','bnb.maxiter', bnbmaxiter);

% set up constraints
Constraints=[...
   ...
   [[W,    H'    ];
    [H  , eye(n,n)]
   ]>=0, ...
   ...
   mu>=0, eta>=0, ...
   ...
   [[eye(n,n) – mu*inv(A) ,                H*inv(B)                    ];
    [inv(B)*H'             , eye(m,m) - mu*inv(B) + inv(B)*W*inv(B)   ]
   ] >= 0, ...
   ...
   [[eye(n,n) + eta*inv(A),                H*inv(B)                    ];
     [inv(B)*H'            , eye(m,m) + eta*inv(B) + inv(B)*W*inv(B)  ]
   ] >= 0, ...
   ...
   H==DA*Htilde*DB, sum(Htilde(:))>=1, ...
   ...
   H*[zeros(kB,m-kB); eye(m-kB,m-kB)]==zeros(n,m-kB), ...
   ...
   Haux==DA*H*DB, sum(Haux(:,:))==diag(W)', ...
   ...
   sum(inv(D)*[[A, H]; [H', B]]*inv(D))<=maxdegree*ones(1,n+m), ...
   ];

% run Yalmip solver
solvesdp(Constraints, -mu-eta, options)

fprintf('\nHOMO-LUMO gap = %f\n', double(mu+eta) );
```

## 5. Computational results

In this section we present various computational examples. We compare results of optimal bridging for weighted and non-weighted graphs representing the structure of chemical molecules of fulvene and benzene. Notice that the weights should correspond to the Coloumb and resonance integrals in real chemical applications (Straitwieser, 1961). For the sake of simplicity, we chose artificial weights illustrating different behavior of optimal bridging for weighted and non-weighted graphs.

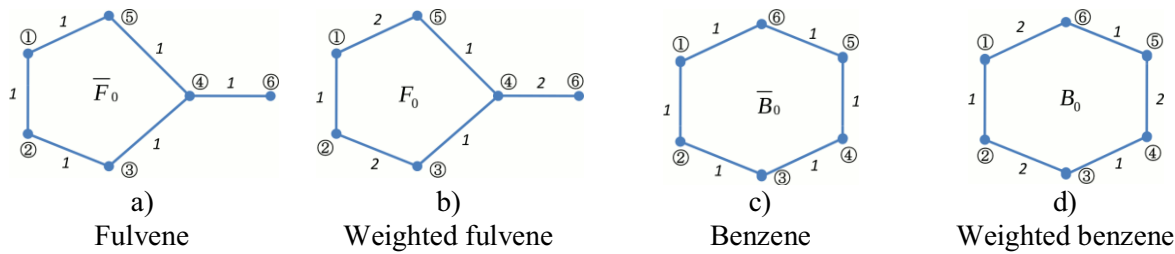

a) Fulvene   b) Weighted fulvene   c) Benzene   d) Weighted benzene

**Figure 2** Structural vertex labeled graphs and weights of edges.

In figure 2, a) we show the graph $\overline{F}_0$ on 6 vertices representing the fulvene organic molecule (5-methylidenecyclopenta-1,3-diene) without weights. The graph $F_0$ with weights is shown in figure 2, b). The voltage matrix for the graph $F_0$ has the form: $D^A = diag(1,1,2,0.5,2,4)$. The structural graph $\overline{B}_0$ of benzene without weights is shown in figure 2, c) and the one with weights $B_0$ is shown in figure 2, d). The voltage matrix for the graph $B_0$ has the form: $D^A = diag(1,1,2,1,1,2)$. Their spectrum consists of the following eigenvalues:

$\sigma(\overline{F}_0) = \{2.1149, 1, 1/q, -0.2541, -q, -1.8608\}$, where $q = (\sqrt{5}+1)/2$,
$\sigma(F_0) = \{3.0680, 1.7437, 1.5616, -1.105, -2.5616, -2.7067\}$,
$\sigma(\overline{B}_0) = \{2, 1, 1, -1, -1, -2\}$,
$\sigma(B_0) = \{3.3723, 2.3723, 1, -1, -2.3723, -3.3723\}$.

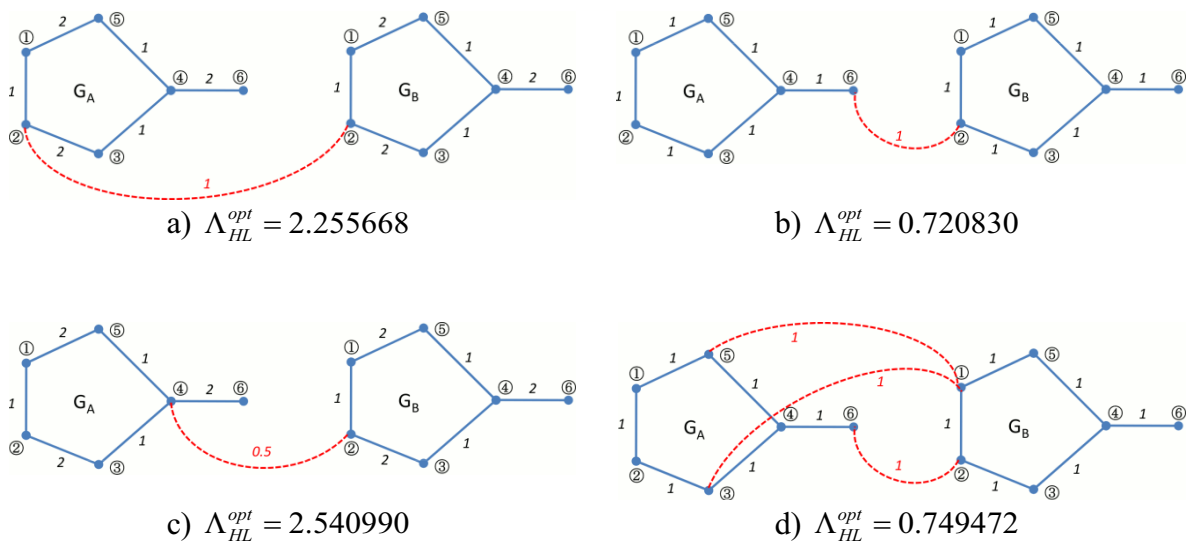

a) $\Lambda_{HL}^{opt} = 2.255668$   b) $\Lambda_{HL}^{opt} = 0.720830$

c) $\Lambda_{HL}^{opt} = 2.540990$   d) $\Lambda_{HL}^{opt} = 0.749472$

**Figure 3** Results of optimal bridging of the graph $G_B$ through the vertices $\{1,2\}$ to the graph $G_A$. The columns a), c) contain graphs $G_A = G_B = F_0$ where $F_0$ is the weighted voltage graph of fulvene; the columns b), d) contain graphs $G_A = G_B = \overline{F}_0$ where $\overline{F}_0$ is the non-weighted fulvene graph. The rows a), b) illustrate optimal bridging with an additional constraint on the vertex degree less or equal to 3; the rows c), d) contain the optimal bridged graphs without constraint on vertex degrees.

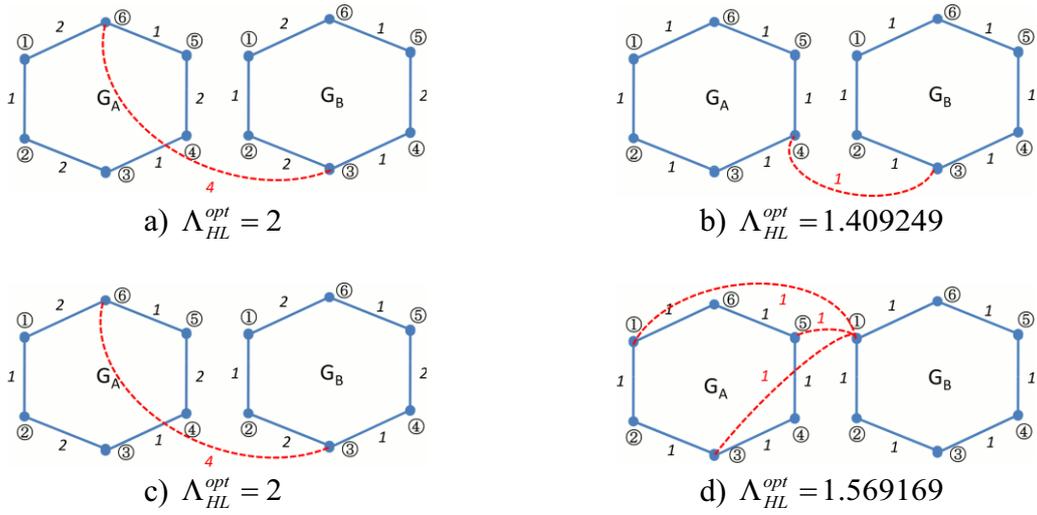

**Figure 4** Results of optimal bridging of the graph $G_B$ through the vertices $\{1,3\}$ to the graph $G_A$. The columns a), c) contain graphs $G_A = G_B = B_0$ representing the weighted voltage graph $B_0$ of benzene; the columns b), d) contain graphs $G_A = G_B = \overline{B_0}$ representing the non-weighted graph $\overline{B_0}$. The rows a), b) illustrate optimal bridging with an additional constraint on the vertex degree less or equal to 3; the rows c), d) contains the optimal bridged graphs without any constraint on vertex degrees.

Recall that a structured graph represents the so-called chemical molecule if its adjacency matrix contains only vertices of the degree less or equal to 3. Such a constraint can be easily implemented into the mixed integer semidefinite programming problem (14) by adding the inequality:

$$\sum_{\substack{j=1\\j\neq i}}^{n+m} \overline{C}_{ij} \leq 3 \quad \text{where} \quad \overline{C} = \begin{pmatrix} \overline{A} & \overline{H} \\ \overline{H}^T & \overline{B} \end{pmatrix}.$$

In table 2 below we present results of optimal bridging between weighted and non-weighted graphs. The results were computed on a Quad Core Intel 1.5GHz CPU with 4GB of memory.

**Table 2** The computational results and comparison of optimal bridging of weighted and non-weighted fulvene and benzene graphs. The first two columns describe the graph $G_A$ and $G_B$ with the chosen set of bridging vertices. The next column indicates whether the constraint on the maximal degree has been imposed. The optimal value $\Lambda_{HL}^{opt}$ is shown in the fourth column.

| $G_A$ | $G_B$ | Max vertex degree | $\Lambda_{HL}^{opt}$ | Comp. time (sec) | Optimal bridging $G_A \to G_B$ and (weights) |
|---|---|---|---|---|---|
| $F_0$ | $F_0$ | 3 | 2.255668 | 10 | $1 \to \emptyset;\ 2 \to 2(1)$ |
| $\bar{F}_0$ | $F_0$ | 3 | 0.835521 | 8.4 | $1 \to \emptyset;\ 2 \to 2(1)$ |
| $F_0$ | $\bar{F}_0$ | 3 | 0.835521 | 21.9 | $1 \to \emptyset;\ 2 \to 2(1)$ |
| $\bar{F}_0$ | $\bar{F}_0$ | 3 | 0.720830 | 9.2 | $1 \to \emptyset;\ 2 \to 6(1)$ |
| $F_0$ | $F_0$ | --- | 2.540990 | 81.2 | $1 \to \emptyset;\ 2 \to 4(0.5)$ |
| $\bar{F}_0$ | $F_0$ | --- | 0.871933 | 29.2 | $1 \to 3(1), 5(1), 6(1); 2 \to \emptyset$ |
| $F_0$ | $\bar{F}_0$ | --- | 1.140215 | 578 | $1 \to 1(1), 2(1), 5(2);$ $2 \to 1(1), 2(1), 3(2)$ |
| $\bar{F}_0$ | $\bar{F}_0$ | --- | 0.749472 | 238 | $1 \to 3(1), 5(1);\ 2 \to 6(1)$ |
| $B_0$ | $B_0$ | 3 | 2 | 10.2 | $1 \to \emptyset;\ 3 \to 6(4)$ |
| $\bar{B}_0$ | $B_0$ | 3 | 1.666731 | 15.2 | $1 \to \emptyset;\ 3 \to 6(2)$ |
| $B_0$ | $\bar{B}_0$ | 3 | 1.666731 | 19.2 | $1 \to \emptyset;\ 3 \to 3(2)$ |
| $\bar{B}_0$ | $\bar{B}_0$ | 3 | 1.409249 | 24.9 | $1 \to \emptyset;\ 3 \to 4(1)$ |
| $B_0$ | $B_0$ | --- | 2 | 4.7 | $1 \to \emptyset;\ 3 \to 6(4)$ |
| $\bar{B}_0$ | $B_0$ | --- | 2 | 10.2 | $1 \to \emptyset; 3 \to 1(2), 3(2), 5(2)$ |
| $B_0$ | $\bar{B}_0$ | --- | 1.830242 | 164 | $1 \to \emptyset;\ 3 \to 1(1), 5(1)$ |
| $\bar{B}_0$ | $\bar{B}_0$ | --- | 1.569169 | 129 | $1 \to 1(1), 3(1), 5(1); 3 \to \emptyset$ |

## 6. Conclusions

We analyzed spectral properties of voltage graphs which are constructed from two given voltage invertible graphs by bridging them over a bipartite graph. We showed that the HOMO-LUMO spectral gap can be computed by means of a solution to the mixed integer semidefinite programming problem. We investigated the optimization problem in which we constructed the optimal bridging graph by means of maximizing the HOMO-LUMO spectral gap. Several computational examples were presented in this paper.

**Acknowledgements.** The research was supported by the APVV Research Grant 0136-12 (SP) and VEGA grant 1/0780/15 (DŠ).

**References.**